\documentclass[12pt, a4paper]{amsart}
\usepackage{amssymb,amscd}
\usepackage[bitstream-charter]{mathdesign}
\DeclareSymbolFont{cmsymbols}{OMS}{cmsy}{m}{n}
\SetSymbolFont{cmsymbols}{bold}{OMS}{cmsy}{b}{n}
\DeclareSymbolFontAlphabet{\mathcal}{cmsymbols}

\usepackage[
unicode=true,
plainpages = false,
pdfpagelabels,
bookmarks=true,
bookmarksnumbered=true,
bookmarksopen=true,
breaklinks=true,
backref=false,
colorlinks=true,
linkcolor = blue,
urlcolor  = blue,
citecolor = red,
anchorcolor = green,
hyperindex = true,
hyperfigures
]{hyperref}
\hypersetup{
pdftitle={Linear algebra, theory and algorithms},
pdfauthor={V. Mikaelian},
pdfsubject={Linear Algebra}
}

\newcommand{\V}{{\mathfrak V}}

\newcommand{\U}{{\mathfrak U}}
\newcommand{\A}{{\mathfrak A}}
\newcommand{\Wr}{\,\mathrm{Wr}\,}

\newcommand{\var}[1]{\mathrm{var}\left( #1 \right)}
\newcommand{\varr}[1]{\mathrm{var}( #1 )}

\theoremstyle{plain}
\newtheorem{Theorem}{Theorem}[section]

\theoremstyle{definition}

\theoremstyle{remark}
 
\newtheorem{Example}[Theorem]{Example} 

\numberwithin{equation}{section}

\begin{document}

\title{Comparing varieties generated by certain wreath products of groups}

\author[V.\,H. Mikaelian]{V.\,H. Mikaelian}

\maketitle

{
\hrule
\vskip2mm
\footnotesize
\noindent
An exhaustive version of the thesis to a talk presented at the: \\
{Groups \& Algebras in Bicocca Conference (GABY)} \\
University of Milano-Bicocca (Milan, Italy),
June 17 to June 21, 2024.\\
\href{https://staff.matapp.unimib.it/~/gaby/gaby2024/index.html}{\tt \footnotesize https://staff.matapp.unimib.it/\textasciitilde/gaby/gaby2024/index.html}
\vskip2mm
\hrule
}

\bigskip
\section{Introduction}
\noindent
In the current talk we would like to  present the main results of recent work \cite{M2018} in which we  classify cases when \textit{the wreath products of distinct pairs of groups  generate the same variety}. This classification allows us to study the subvarieties of some nilpotent-by-abelian product varieties  $\U\V$ with the help of wreath products of groups.
For  background information on varieties of groups, on generating groups of varieties, on products of varieties, or on wreath products we refer to Hanna Neumann's 
monograph \cite{HannaNeumann}, to the related articles
\cite{3N,
Some_remarks_on_varieties,
Baumslag nilp wr,
B+3N, 
ShmelkinOnCrossVarieties, Burns65},
and to literature cited therein.

In particular,  using wreath products we discover such subvarieties in nilpotent-by-abelian products $\U\V$ which \textit{have the same nilpotency class, the same length of solubility, and the same exponent, but which still are distinct subvarieties}. Obtained classification strengthens results on varieties generated by wreath products in the mentioned articles and elsewhere in the literature, see \cite{M2018} for references.

\medskip
Wreath products are among the main tools to study products $\U\V$ of any varieties $\U$ and $\V$ of groups. 
Under wreath products we by default suppose \textit{Cartesian} wreath products, but all the results we bring are true for \textit{direct} wreath products also. 
In the literature the wreath product methods most typically consider certain groups $A$ and $B$ generating the varieties $\U$ and $\V$ respectively, and then they find extra conditions, under which the wreath product ${A \Wr B}$ generates $\U\V$, i.e., conditions, under which the equality
\begin{equation}
\tag{$*$} 
\label{EQUATION_main}    
\var{A \Wr B} = \var{A} \var{B}
\end{equation}
holds for the given $A$ and $B$. 

The advantage of such an approach is that having the equality \eqref{EQUATION_main} we using Birkhoff's Theorem can get \textit{all} the groups in $\var{A} \var{B}=\U\V$ by applying the operations of taking the homomorphic images, subgroups, Cartesian products to the \textit{single} group $A \Wr B$ only, see  \cite[15.23]{HannaNeumann}.

Generalizing some known results in the cited literature, we in \cite{AwrB_paper}, 
\cite{M2007-b},
\cite{Algebra i logika - first},
\cite{M2017},
\cite{K_p-series},
\cite{Classification Theorem} were able to suggest criteria classifying all the cases when \eqref{EQUATION_main} holds for groups from certain particular classes of groups: abelian groups, $p$-groups, nilpotent groups of finite exponent, etc. See the very brief outline of results presented in Section 5 of \cite{Classification Theorem}.

\medskip
In \cite{M2018} we turned to a sharper problem of comparison of two varieties, both generated by wreath products
. Namely, take $A_1, B_1$ and $A_2, B_2$ to be pairs of non-trivial groups such that\,
$\varr {A_1}=\varr {A_2}$, \, $\varr {B_1}=\varr {B_2}$,\, and then dis\-tinguish the cases when:
\begin{equation}
\tag{$**$}
\label{Equation varieties are equal}
\var{A_1 \!\Wr B_1} = \var{A_2 \!\Wr B_2}.
\end{equation}

\bigskip
\section{The main criterion and examples}
\noindent
To write down the main classification criterion 
from \cite{M2018} we need some very simple notation.
Namely, by Pr{\"u}per's Theorem any abelian group $B$ of finite exponent is a direct product of its $p$-primary components $B(p)$, and each of such components is a direct product of certain cyclic $p$-groups $C_{p^{u_1}}, C_{p^{u_2}},\ldots$ 
If among the latters the cardinality of copies isomorphic to $C_{p^{u_k}}$ is $m_{p^{u_k}}$, we can write their direct product as  $C_{p^{u}}^{m_{p^{u}}}$\!. Then the component $B(p)$ can be rewritten as a direct product of such factors:
\begin{equation}
\label{Equation decomposition of B(p)}
B(p)
=C_{p^{u_1}}^{m_{p^{u_1}}}
\!\!\times \cdots \times
C_{p^{u_r}}^{m_{p^{u_r}}}\!\!\!,
\end{equation}
where we may suppose $u_1 > \cdots > u_r$, see \cite[Section 35]{Fuchs1}. 
If $B(p)$ is finite, then all the cardinals $m_{p^{u_1}},\ldots,m_{p^{u_r}}$ together with all the factors $C_{p^{u_k}}^{m_{p^{u_k}}}$ are finite. Otherwise, at least one of those factors has to be infinite, and we can denote $C_{p^{u_k}}^{m_{p^{u_k}}}$
to be the \textit{first one} of them, i.e., 
$m_{p^{u_k}}$ is an infinite (at least countable) cardinal, and all the preceding cardinals $m_{p^{u_1}},\ldots,m_{p^{u_{k-1}}}$ are finite.

Let $B_1$ and $B_2$ be abelian groups of finite exponents divisible by some prime $p$.
Call their $p$-primary components $B_1(p)$ and $B_2(p)$ 
\textit{equivalent}, if in \eqref{Equation decomposition of B(p)} 
their first \textit{infinite} direct factors have the same exponent, and all their preceding \textit{finite} direct factors coincide. 
More precisely, if 
$B_1(p)
=C_{p^{u_1}}^{m_{p^{u_1}}}
\!\!\times \cdots \times
C_{p^{u_r}}^{m_{p^{u_r}}}$
and
$B_2(p)
=C_{p^{v_1}}^{ m_{p^{v_1}}}
\!\!\times \cdots \times
C_{p^{v_s}}^{m_{p^{v_s}}}$\!, then $B_1(p)\equiv B_2(p)$ if and only if:
\begin{enumerate}
\item when $B_1(p),B_2(p)$ are \textit{finite}, then  $B_1(p) \equiv B_2(p)$ iff $B_1(p) \cong B_2(p)$;
\item when $B_1(p),B_2(p)$ are \textit{infinite}, then  $B_1(p) \equiv B_2(p)$ iff there is a $k$ so that:
\begin{enumerate}
\item[\textit{i}) ] $u_i=v_i$ and ${m_{p^{u_i}}}={m_{p^{v_i}}}$ for each $i=1,\ldots,k-\!1$;
\item[\textit{ii}) ] $C_{p^{u_k}}^{ m_{p^{u_k}}}$\!\! is the first infinite factor for $B_1(p)$;\; $C_{p^{v_k}}^{m_{p^{v_k}}}$\!\! is the first in\-finite factor for $B_2(p)$, and $u_k=v_k$;
\end{enumerate}
\item else $B_1(p),B_2(p)$ are \textit{not} equivalent.
\end{enumerate}
The above definition is not short, but it is very intuitive to understand:

\begin{Example}
\label{EX one}
$C_{3^5}^{6} \times 
C_{3^4}^{8} \times 
C_{3^3}^{\aleph_0}\times 
C_{3^2}^{5} \times 
C_{3}^{4}$
is equivalent to 
$C_{3^5}^{6} \times 
C_{3^4}^{8} \times 
C_{3^3}^{\mathfrak c}\times 
C_{3}^{50}$, but it is \textit{not} equivalent to 
$C_{3^5}^{6} \times 
C_{3^4}^{8} \times 
C_{3^2}^{\aleph_0} \times 
C_{3}^{4}$. Here 
$\aleph_0$ and $\mathfrak c$ stand for countable and continuum cardinals.
In the first two of the above groups the first \textit{infinite} factors  
$C_{3^3}^{\aleph_0}$ and $C_{3^3}^{\mathfrak c}$ are of the same exponent $3^3$ (without being isomorphic), and they both have the same two initial \textit{finite} factors $C_{3^5}^{6} \times 
C_{3^4}^{8}$.
Whereas in the third group the first \textit{infinite} factor $C_{3^2}^{\aleph_0}$ is of another exponent $3^2 \neq 3^3$.
\end{Example}

In these terms our main criterion reads:
\begin{Theorem}
\label{Theorem var(Wr) are equal}
Let $A_1, A_2$ be any non-trivial nilpotent groups of  exponent $m$ generating the same variety, and let $B_1,\,B_2$ be any non-trivial abelian groups of exponent $n$  generating the same variety, where any prime divisor $p$ of $n$ also divides $m$. 

Then \eqref{Equation varieties are equal} holds for $A_1, A_2, B_1, B_2$ if and only if $B_1(p) \equiv B_2(p)$ for each $p$.
\end{Theorem}

Notice how the roles of the passive and active groups of these wreath products are different: for $A_1, A_2$ we just require that $\var{A_1} = \var{A_2}$, whereas for $B_1, B_2$ we put extra conditions on
structures of their decompositions. And when 
$B_1, B_2$ are finite, the extra conditions simply mean $B_1 \cong B_2$.

\begin{Example}
\label{EX two}
To see an application of Theorem~\ref{Theorem var(Wr) are equal} take $Q_8$ to be the quaternion group, 
and take $M_{27}$ to be the semidirect product of $C_9$ and of $C_3$, acting on it by nontrivial automorphisms.
$Q_8$ is of order $8$, of exponent $4$, and of nilpotency class $2$, while
$M_{27}$ is of order $27$, of exponent $9$, and of nilpotency class $2$.
Then pick 
$A_1= Q_8 \times M_{27} \times C_{25}$,\;\;
$A_2= Q_8 \times Q_8 \times M_{27}^{\aleph_0} \times C_{25}\times C_5\times C_5$,\;\;
$B_1=C_{3^5}^{6} \times 
C_{3^4}^{8} \times 
C_{3^3}^{\aleph_0}\times 
C_{3^2}^{5} \times 
C_{3}^{4} \times C_5$,\;\;
$B_2=C_{3^5}^{6} \times 
C_{3^4}^{8} \times 
C_{3^3}^{\mathfrak c}\times 
C_{3}^{50} \times C_5$.
Three conditions 
$\var{A_1}=\var{A_2}$,\;
$B_1(3)\equiv B_2(3)$,\;
$B_1(5)\equiv B_2(5)$
are easy to verify, see Example~\ref{EX one}. Hence \eqref{Equation varieties are equal} holds for this choice of $A_1, A_2, B_1, B_2$.\;
On the other hand, we will no longer have an equality choosing either
$B_2=C_{3^5}^{6} \times 
C_{3^4}^{8} \times 
C_{3^2}^{\aleph_0} \times 
C_{3}^{4} \times C_5$ (because
$3^3 \neq 3^2$), or $B_2=C_{3^5}^{6} \times 
C_{3^4}^{8} \times 
C_{3^3}^{\mathfrak c}\times 
C_{3}^{50} \times C_5 \times C_5$ (because
$(C_5\! \times \!C_5) \not\cong C_5$).
\end{Example}

\bigskip
\section{Applications to subvariety structures}
\noindent
Theorem~\ref{Theorem var(Wr) are equal} 
covers the cases of nilpotent $A_1,A_2$ and abelian $B_1,B_2$, with some restrictions on exponents.
%
Besides getting a generalization of \eqref{EQUATION_main} our study of equality \eqref{Equation varieties are equal} is motivated by some applications one of which we would like to outline here.

Classification of subvariety structures of $\U\V$ is incomplete even when $\U$ and $\V$ are such ``small'' varieties as the abelian varieties $\A_m$ and $\A_n$ respectively. Here are some of the  results in this direction:
$\A_p$ (for prime numbers $p$) are the simplest non-trivial varieties,  as they consist of the Cartesian powers of the cycle $C_p$ only. L.G.~Kov\'acs and M.F.~Newman in~\cite{KovacsAndNewmanOnNon-Cross} fully described the subvariety structure in the product $\A_p^2 = \A_p \A_p$ for $p>2$. 
Later they continued this classification for the varieties $ \A_{p^u} \A_p$. Their research was unpublished for many years, and it appeared in 1994 only \cite{Kovacs Newman Ravello} (parts of their proof are present in \cite{Bryce_Metabelian_groups}).
Another direction is  description of subvarieties in the product $\A_m \A_n$ where $m$ and $n$ are \textit{coprime}. This is done by C.~Houghton (mentioned by Hanna Neumann in \cite[54.42]{HannaNeumann}), by
P.\,J.~Cossey
(Ph.D. thesis \cite{Cossey66}, mentioned by R.A.~Bryce in \cite{Bryce_Metabelian_groups}).
A more general result of R.A.~Bryce classifies the subvarieties of $\A_m \A_n$, where $m$ and $n$ are \textit{nearly prime} in the sense that, if a prime $p$ divides $m$, then $p^2$ does not divide $n$ \cite{Bryce_Metabelian_groups}.

In 1967 Hanna Neumann wrote that classification of subvarieties of  $\A_m \A_n$ for arbitrary $m$ and $n$  \textit{``seems within
reach''} \cite{HannaNeumann}. 
And R.A.~Bryce in 1970 mentioned that \textit{``classifying all metabelian varieties
is at present slight'' } \cite{Bryce_Metabelian_groups}. However, nearly  half a century later this task is not yet accomplished: Yu.A.~Bakhturin and A.Yu.~Olshanskii remarked in the  survey \cite{Bakhturin Olshanskii} of 1988 (appeared in English in 1991) that \textit{``classification of all nilpotent metabelian group varieties has not been completed yet''.}

\smallskip
As this brief summary shows, one of the cases, when the subvariety structure of $\U\V$ is \textit{less known}, is the case when $\U$ and $\V$ have non-coprime exponents divisible by \textit{high} powers $p^u$ for \textit{many} prime numbers $p$.
Thus, even if we cannot classify all the subvarieties in some product varieties $\U\V$, it may be  interesting to find those subvarieties in $\U\V$, which are generated by wreath products. 
We, surely, can take any groups $A \in \U$ and $B \in \V$, and then $A \Wr B$ will generate some subvariety in $\U\V$. 
But in order to make this approach reasonable, we yet have to detect \textit{if or not two wreath products of that type generate the same subvariety}, i.e, if or not the equality \eqref{Equation varieties are equal} holds for the given pairs of groups.

\smallskip
Yet another outcome of this research may be stressed.
In the literature the different subvarieties are often distinguished by their different nilpotency classes, different lengths of solubility, or different exponents (see, for example, classification of subvarieties of $\A_p^2$ in~\cite{KovacsAndNewmanOnNon-Cross}).
Using wreath products technique, we in \cite{M2018} construct such subvarieties of $\U\V$, which have the same nilpotency class, the same length of solubility, the same exponent, but which still are distinct subvarieties, see examples in \cite{M2018}.
Other related research can be foind in 
\cite{M2002-b},
\cite{M2003},
\cite{M2004},
\cite{M2005},
\cite{M2023}.

\medskip
\noindent 
\footnotesize
E-mail:
\href{mailto:v.mikaelian@gmail.com}{v.mikaelian@gmail.com}, 
\href{mailto:vmikaelian@ysu.am}{vmikaelian@ysu.am}
$\vphantom{b^{b^{b^{b^b}}}}$

\noindent 
Web: 
\href{https://www.researchgate.net/profile/Vahagn-Mikaelian}{researchgate.net/profile/Vahagn-Mikaelian}

\end{document}